\documentclass[10pt]{article}

\usepackage{fancyhdr}
\usepackage{extramarks}
\usepackage{amsmath}
\usepackage{amsthm}
\usepackage{amsfonts}
\usepackage{siunitx}
\usepackage{tikz}
\usepackage[plain]{algorithm}
\usepackage{algpseudocode}
\usepackage{multirow}
\usepackage{booktabs}
\usepackage{graphicx}
\usepackage{subfigure}
\usepackage[colorlinks,linkcolor=black,anchorcolor=black,citecolor=black,urlcolor=blue]{hyperref}
\usepackage{amsmath,bm}
\usepackage{booktabs}
\usepackage{mathtools}
\usepackage{amssymb}
\usepackage{caption}
\usepackage{capt-of}
\usepackage{mciteplus}
\usepackage{cite}
\usepackage{mathrsfs}
\usepackage[title,titletoc,toc]{appendix}
\usepackage{xr}
\usepackage{parskip}
\usetikzlibrary{automata,positioning}
\usepackage{wrapfig}
\renewcommand{\part}[1]{\textbf{\large Part \Alph{partCounter}}\stepcounter{partCounter}\\}

\begin{document}

\title{Topological data analysis ``hearing’’ the shapes of drums and bells   }
\author{By  Guo-Wei Wei 
}
%\date{\today} % Date for the report

\maketitle
%\linenumbers 
\begin{abstract} 

Mark Kac asked a famous question in 1966, ``Can one hear the shape of a drum?’’, a spectral geometry problem that has intrigued mathematicians for the last six decades and is important to many other fields, such as architectural acoustics, audio forensics, pattern recognition, radiology, and imaging science.  A related question is how to hear the shape of a drum. We show that the answer was given in the set of 65 Zenghouyi chime bells dated back to 475-433 B.C. in China. The set of chime bells gradually varies their sizes and weights to enable melodies, intervals, and temperaments. The same design principle was used in many other musical instruments, such as xylophones, pan flutes, pianos, etc. We reveal that  there is a fascinating connection between the progression pattern of many musical instruments and filtration (or spectral sequence) in topological data analysis (TDA). We argue that filtration-induced evolutionary de Rham-Hodge theory provides a new mathematical foundation for musical instruments. Its discrete counterpart, persistent Laplacians and many other persistent topological Laplacians, including persistent sheaf Laplacians  and persistent path Laplacians are briefly discussed.   

\end{abstract}

\begin{wrapfigure}{r}[0.0in]{3in} 
	\vspace{-0mm}
	\includegraphics[keepaspectratio,width=3in]{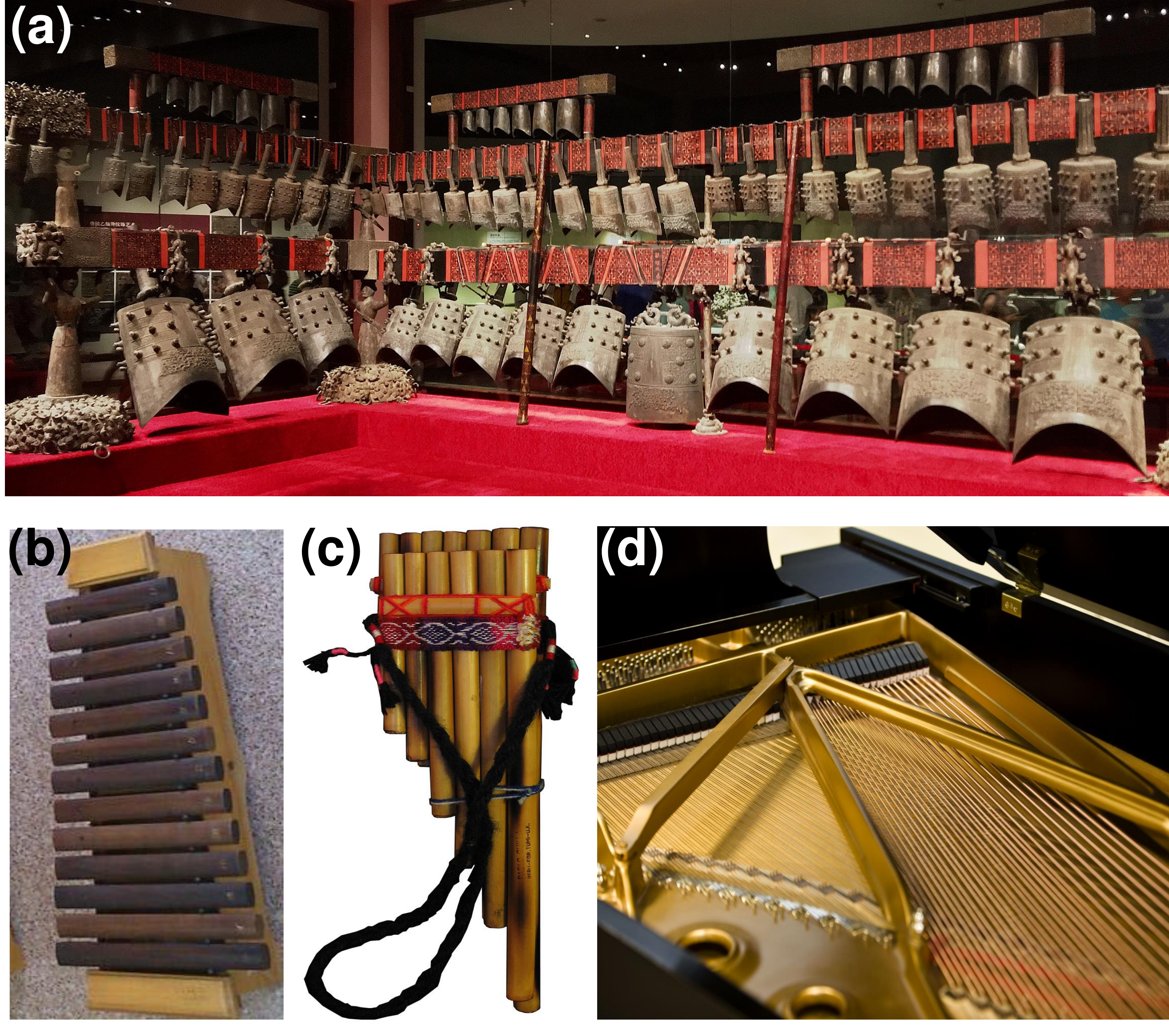}
	\vspace{-5mm}
    \caption{The progression patterns in musical instruments. 
			{\bf a}. The Zonghouyi chime bells. 
			{\bf b}. An xylophone.
			{\bf c}. A pan flute.
			{\bf d}. Strings of a grand piano.   
			  Image courtesy of Wikipedia.  
		}
    \label{fig:music}
\end{wrapfigure}
 ``Can one hear the shape of a drum?’’,  a famous question posed by Mark Kac in 1966 \cite{kac1966can}, has intrigued mathematicians for generations. In other words, if you hear the sound from a drum, i.e., the set of overtones produced by the drum, can you infer its shape?   
Mathematically, the essence of the spectral geometry question is whether the shape can be uniquely determined from the eigenvalues of the Laplacian operator defined on the shape. There are many examples of isospectral manifolds which are not isometric in   two    and higher dimensional settings \cite{milnor1964eigenvalues,gordon1992isospectral, vigneras1980varietes}. However, the problem is not closed yet ---  one can discern the shapes of certain geometric types from their sounds \cite{reid1992isospectrality,hezari2021dirichlet}. The question has a far-reaching impact on many fields beyond mathematics, such as architectural acoustics, audio forensics,  pattern recognition, radiology, imaging science, and musical science \cite{dokmanic2013acoustic,cosmo2019isospectralization}. 

In musical composition, it is a common practice to use a drum set of varying sizes, instead of a single piece of drum, to facilitate a tonic harmonic progression, which is a foundation of harmony in modern music. The human brain is trained from the drum set to distinguish the overtones of individual drum. The comparative training/learning from a set of drums is the basis to hear the shape of a drum from a drum set. 

Interestingly, to create tonal harmony, the ancient Chinese built a set of 65 Zenghouyi chime bells dated back to 475-433 B.C.  in the Warring States Period in China (Figure 1a). The shape of each chime bell was designed to produce a distinct sound.  The set of 65 chime bells gradually varies in their sizes and weights, ranging from 153.4 centimeters (60.4 in) to 20.4 centimeters (8.0 in) in height and from 203.6 kilograms (449 lb) to 2.4 kilograms (5.3 lb) in weight. The set of Zenghouyi chime bells covers from C2 to D7 tonal range and can play all twelve half tones in the middle area of the tonal range \cite{yan2013physics}, enabling melodies, intervals, and temperaments. The same design principle is used in many other musical instruments, such as xylophones (Figure 1b), pan flutes (Figure 1c), and  pianos (Figure 1d).  
Due to the close proximity among chime bells, resonance among them may occur  when they are struck, giving rise to prolonged harmony. Made of bronze, one of the most precious metals available at the time, chime bells were used in various rituals, ceremonies, and  entertainment. 

There is a fascinating and apparent connection between the progression pattern of musical instruments (Figure 1) and mathematical filtration (Figure 2a), a 
spectral sequence technique \cite{mccleary2001user} widely used in homological algebra and topological data analysis (TDA). As a new branch of mathematics,  TDA uses topological and geometric concepts to understand and extract topological patterns and structures in data. The main workhorse of TDA is persistent homology \cite{edelsbrunner2008persistent,Zomorodian:2005}, a branch of algebraic topology that creates a mathematical microscopy of a point cloud by filtration. Through the comparative analysis of the topological invariant changes induced by the filtration, persistent homology delineates the shape of data \cite{lum2013extracting}. Paired with advanced machine learning and deep learning algorithms, persistent homology has had tremendous success in data science.  It has been established as one of the most powerful tools in simplifying the geometric complexity and reducing the high dimensionality of biomolecular interactions \cite{liu2022biomolecular}, revolutionizing \href{https://sinews.siam.org/Details-Page/persistent-homology-analysis-of-biomolecular-data}{drug discovery} and  \href{https://sinews.siam.org/Details-Page/topological-artificial-intelligence-forecasting-of-future-dominant-viral-variants}{the forecasting of emerging viral variants}.  
  
\begin{wrapfigure}{r}[0.0in]{3.in} 
	\vspace{-2mm}
	\includegraphics[keepaspectratio,width=3.in]{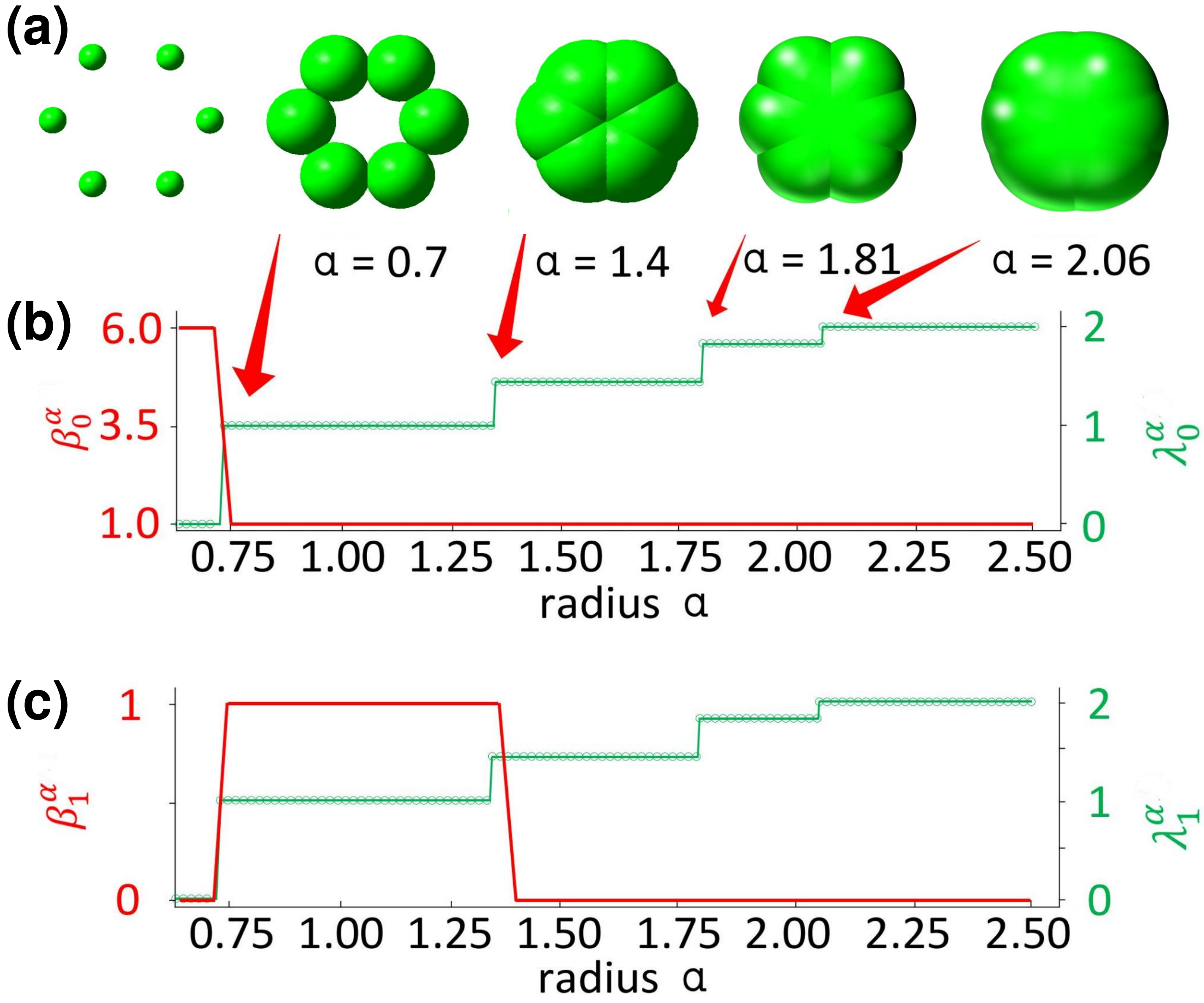}
	\vspace{-5mm}
    \caption{Illustration of  filtration and persistent Laplacian spectra. 
			{\bf a}.  Filtration of a point cloud. 
			{\bf b}-{\bf c}. The persistent Laplacian analysis  of the point cloud. 
			 Here, $\beta^{\alpha}_j$ and $\lambda^{\alpha}_j$ ($j=0,1$) are  persistent Betti-$j$ numbers and the first nonzero eigenvalues of the $j$th persistent Laplacian, respectively.  Note that $\lambda^{\alpha}_j$ capture the homotopic shape evolution of data (see the frequency changes after $\alpha=1.5$) that is not reflected in $\beta^{\alpha}_j$. The similarity between $\lambda^{\alpha}_0$ and $\lambda^{\alpha}_1$ is a coincidence. 
			 Image courtesy of Dr. Jian Jiang.  
		}
    \label{fig:Laplacian}
\end{wrapfigure}
However, persistent homology is not directly applicable to the tonal analysis of chime bells. First, the set of chime bells may be regarded as the result of a set of evolving chime bell manifolds, rather than that of the filtration of a point cloud. Additionally, there is no change in topological invariants associated the homotopic shape evolution of the set of chime bells. Finally, persistent homology cannot present a frequency or tonal analysis of chime bells or many other musical instruments.    
  
Recently, an evolutionary de Rham-Hodge method has been proposed as a multiscale generalization of the classical de Rham-Hodge theory, a landmark of the 20th Century’s mathematics \cite{chen2021evolutionary}. It provides a multiscale geometric and topological analysis of filtration-induced manifolds, on which a family of evolutionary Hodge Laplacians can be defined on the set of chime bells to characterize their tonal evolution.  In association with a family of evolutionary de Rham complexes,  evolutionary Hodge Laplacians reveal the full set of topological invariants, such as Betti numbers, in their kernel or null space dimensions. 

The point-cloud counterpart of the evolutionary de Rham-Hodge method on manifolds is called persistent spectral graph \cite{wang2020persistent} (also known as persistent Laplacians \cite{memoli2022persistent}) on simplicial complexes. Like the evolutionary de Rham-Hodge method, persistent Laplacians not only return  the full set of topological invariants in their harmonic spectra as persistent homology does but also capture the homotopic shape evolution of data during the filtration in their first non-harmonic spectra  (Figure 2b and 2c), for which persistent homology cannot describe.    

A generalization of persistent Laplacians was made through the sheaf theory  \cite{hansen2019toward} and the resulting persistent sheaf Laplacians \cite{wei2021persistent} allow the embedding of heterogeneous characters in topological invariants, e.g., encoding non-geometric information in a geometry-based simplicial complex.  Another generalization is persistent path Laplacians \cite{wang2022persistent}, built from the path complex and path homology \cite{grigor2012homologies,chowdhury2018persistent}. Persistent path Laplacians are designed  for directed graphs and directed networks. These new persistent topological Laplacians not only lay a mathematical foundation for the tonal analysis in musical science but also significantly extend the applicable domain and power of TDA.

%{\setcounter{tocdepth}{4} \tableofcontents}

 % \setcounter{page}{1}
 %\renewcommand{\thepage}{{\arabic{page}}}

 % \bibliographystyle{abbrv}
 
% \bibliography{refs}

\vspace{0.8cm}
{\it Guo-Wei Wei is an MSU Foundation Professor at Michigan State University. His research concerns the mathematical foundations of biological science and data science.}

\end{document}